\documentclass[12pt,oneside,reqno]{amsart}
\usepackage{amsmath,amsfonts,amssymb,graphicx,color}
\usepackage[initials]{amsrefs}
\usepackage[all]{xy}
\usepackage{geometry}
\theoremstyle{definition}

\DeclareMathOperator*{\VV}{[\,\, ]^{-1}}
\DeclareMathOperator*{\VE}{\varepsilon\theta}

\input epsf
\begin{document}
%%%%%%%%%%%%%%
\title[On the invariant subspace problem in Hilbert spaces]{On the invariant subspace problem\\ in Hilbert spaces}
%%%%%%%%%%%%%%
%%%%%%%%%%%%%%
\author[Enflo]{Per H. Enflo}
\address[Per H. Enflo]{\mbox{}\newline \indent Department of Mathematical Sciences,\newline\indent  Kent State University, \newline\indent Kent, Ohio, 44242, USA.\newline \indent
	\textsc{ and }\newline \indent  \noindent Instituto de Matem\'atica Interdisciplinar (IMI), \newline \indent Facultad de Ciencias Matem\'{a}ticas \newline \indent Plaza de Ciencias 3 \newline \indent Universidad Complutense de Madrid \newline \indent	Madrid, 28040 (Spain)
}
\email{per.h.enflo@gmail.com}
%%%%%%%%%%%%%%
%%%%%%%%%%%%%%
%%%%%%%%%%%%%%
%%%%%%%%%%%%%%
%%%%%%%%%%%%%%
\begin{abstract}
In this paper we show that every bounded linear operator $T$ on a Hilbert space $H$ has a closed non-trivial invariant subspace.
\end{abstract} 
%%%%%%%%%%%%%%
%%%%%%%%%%%%%%
\maketitle
%%%%%%%%%%%%%%
%%%%%%%%%%%%%%

The aim of this paper is to show that every bounded linear operator $T$ on a Hilbert space $H$ has a closed non-trivial invariant subspace. This gives a solution to the long-standing Invariant Subspace Problem.  

We will construct convergent sequences in $H$, whose limits will be non-cyclic vectors. In the construction there are many cases and choices to consider. We have found it appropriate to present a few of the considerations as Lemmas, but most of them as part of a continuing process of reasoning, represented by numbers in brackets.

Instead  of trying to present a detailed construction plan in the beginning of the paper, we shall give a ``Summary of the Construction'' at the end of the paper, which -- with all definitions and considerations in place -- should make our strategy easier to be followed.

Since the spectrum of $T$, $\sigma (T)$, is non-empty and $\lambda I + \mu T$, $\mu \neq 0$, has the same closed invariant subspaces as $T$, we can (without loss of generality) assume that $T$ is one to one, $\mathcal{R}(T)  \neq H$, $\overline{\mathcal{R}(T)}=H$, $0 \in \sigma (T)$. In order to avoid delicate calculations we will assume that $\|T\|_{\text{op}}$ is small (we put $\|T\|_{\text{op}} = 10^{-20}$ below). 

The method to construct invariant subspaces represents a new direction of the author's method of ``extremal vectors'', first presented in \cite{AnsariEnflo}. Solving norm inequalities has been one of the many important tools when dealing with the invariant subspace problem for Banach spaces (see, e.g., \cite{enflo1987, MR4093902, BBEM}).  

In \cite{AnsariEnflo} one uses the $y \in H$ of minimal norm to solve the inequality $\|x_0 - T^n y\| \le \varepsilon$, $\|x_0\| = 1$, $0 < \varepsilon < 1.$ In this paper, for a given $y'$, we use the $\ell '$, 
$$\ell^\prime (T) = \sum_{j \ge 0} a_j T^j$$
of minimal $\|\cdot \|_2$-norm,
$$\left\| \sum_{j \ge 0} a_j T^j\right\|_2 = \left( \sum_{j \ge 0} |a_j|^2\right)^{1/2}$$
to solve the inequality 

\begin{equation}\label{eq01}
	\|x_0 - \ell^\prime (T) y' \| \le \varepsilon
\end{equation}

Throughout this paper we can think of $y'$ having a norm of the order of magnitude $1$, whereas $y = s y'$, $s>1$, has a bigger norm.

With 
$$r(T) = \sum_{j \ge 0} r_j T^j \quad \text{ and } \quad \ell^\prime (T) = \sum_{j\ge 0} a_j T^j$$

we have $\displaystyle \langle r, \ell^\prime \rangle = \sum_{j \ge 0} r_j \overline{a_j}$.

The minimality of $\|\ell^\prime\|_2$ in \eqref{eq01} gives that

\begin{equation}\label{eq02}
\displaystyle \langle r, \ell^\prime \rangle = 0 \text{ implies } \langle r(T) y', x_0 - \ell^\prime (T) y\rangle = 0
\end{equation}

We put $V_{y'} \ell' = \ell'(T) y'$ and

\begin{equation}\label{eq03}
V^{*}_{y'} x = \langle x, y' \rangle I + \langle x, T y'\rangle T + \ldots = \sum_{j \ge 0} \langle x, T^{j} y'\rangle T^j
\end{equation}

Equation \eqref{eq03} gives

\begin{equation}\label{eq04}
\langle	V^{*}_{y'} x, \ell' \rangle =   \sum_{j \ge 0} \langle x, T^{j} y'\rangle \overline{a_j} 
= \sum_{j \ge 0} \langle x, a_j T^{j} y'\rangle = \langle x, V_{y'}\ell' \rangle
\end{equation}

We observe here that $V^{*}_{y'} x$ is an operator on $H$ and $V_{y'} \ell' \in H$. The minimality of $\|\ell'\|_2$ gives
\begin{equation}\label{eq05}
	\langle	 V_{y'} r , x_0 - V_{y'} \ell' \rangle =  C' \langle r, \ell' \rangle .
\end{equation}
\medskip

If $C'=0$, $\langle r(T) y', x_0 - \ell' (T) y' \rangle = 0$ for every $r$, then 
$$\overline{\text{span}} \left\{ T^j y' \, : \, j \ge 0\right\}$$
is a non-trivial closed, invariant subspace. Thus, let us assume $C' \neq 0.$\\
\medskip

Equation \eqref{eq05} gives, with $r = \ell'$, 
\begin{equation}\label{eq06}
\langle V_{y'} \ell', x_0 - V_{y'} \ell' \rangle = C' \langle \ell', \ell'\rangle
\end{equation}

We will now see that $\langle V_{y'} \ell' , x_0 - V_{y'} \ell'\rangle$ is a positive real number. For, if $\langle V_{y'} \ell' , x_0 - V_{y'} \ell'\rangle$ were not a real number, we could decrease $\|x_0 - V_{y'}\ell'\|$ by multiplying $\ell'$ by a complex number of absolute value $1$, contradicting the minimality of $\|\ell'\|_2$. To see this, put $V_{y'} \ell' = \alpha x_0 + z,$ $\langle z, x_0 \rangle = 0.$ 

\begin{equation}\label{eq07}
	\|x_0 - V_{y'}\ell'\|^2 = |1 - \alpha|^2 + \|z\|^2,
\end{equation}
which, for a given absolute value of $\alpha$, is minimal for real $\alpha$.\\
\medskip

\begin{figure}[h!]
	\centering
\fbox{\includegraphics[height=.35\textwidth,keepaspectratio=true]{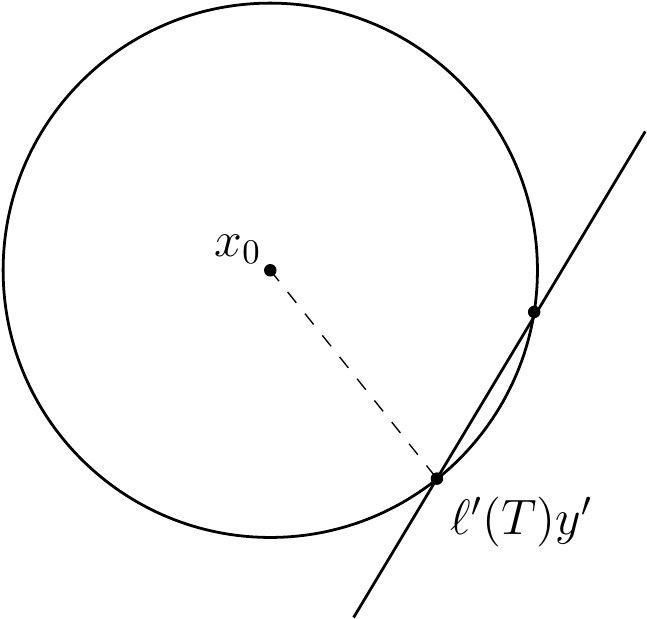}}
\caption{}
\label{figure}
\end{figure}

The minimality of $\|\ell'\|_2$ gives $\langle \ell^\prime (T) y', x_0 -  \ell^\prime (T) y' \rangle \ge 0$. We feel that this can most easily be understood through Figure \ref{figure}.

We put

\begin{equation}\label{eq08}
\langle \ell'(T) y', x_0 - \ell'(T) y'\rangle = \varepsilon \theta.
\end{equation}

With $$T^j \ell' (T) = \alpha_j \ell'(T) + r_j (T), \quad \langle r, \ell' \rangle = 0, \quad |\alpha_j|\le 1$$ we get
\begin{equation}\label{eq09}
	\left| \langle T^j \ell'(T) y', x_0 - \ell'(T) y' \rangle \right|= |\alpha_j| \varepsilon \theta \le \varepsilon \theta.
\end{equation}
for every $j\ge 0$. Thus, if in \eqref{eq08}, $\varepsilon \theta = 0$ then $\ell'(T) y'$ is non-cyclic. Thus, we assume $\varepsilon \theta > 0$.

With  $0.3 \le \| x_0 - \ell' (T) y' \| \le 0.7$, as we will have below, $\|\ell' (T) y'\| \ge 0.3$ and we get for sufficiently small $\delta >0$, 
$$\| x_0 - \ell'(T) y' (1 + \delta) \| \le \|x_0 - \ell' (T) y' \| - \frac{\delta (\varepsilon \theta)}{10},$$
so 
\begin{equation}\label{eq10}
	\text{inf} \left\{ \| \ell\|_2, \, \|\ell(T) y'- x_0\| \right\} \le \|\ell'(T) y' - x_0 \|  - \frac{\delta (\varepsilon \theta)}{10} < \|\ell'\|_2 (1 + \delta).
\end{equation}
In the Main Construction (MC) described below, we have  a sequence $(y_n^\prime)$ such that $0.7 > \varepsilon > 0.3$, $\varepsilon = \|x_0 - \ell'(T) y_n^\prime\|$ for the whole sequence, such that $(\varepsilon \theta)_n \rightarrow 0$ and with $\ell'_n(T) y_n^\prime$ converging  in norm. Then, for the limit $\ell'_{\infty}(T)y_{\infty}^{\prime}$ we will have 
$$\langle T^j \ell'_{\infty}(T) y'_{\infty}, x_0 - \ell'_{\infty}(T) y'_{\infty}\rangle = 0$$
for all $j \ge 0$, i.e., 

\begin{equation}\label{eq11}
	\ell^{\prime}_{\infty}(T) y_{\infty}^{\prime} \text{ is non-cyclic}.
\end{equation}

The approximation of a non-cyclic vector by ``$(\varepsilon \theta)_n$-almost cyclic vectors'', $(\varepsilon \theta)_n \rightarrow 0$, cannot work for certain shift operators. It is well-known that for weighted shifts $e_n \rightarrow w_n e_{n+1}$ , $(e_n)$ ON basis of $H$, with $\displaystyle \sum_n \|w_n\|^2 < \infty$ the only non-cyclic vectors $\sum_{j \ge 0} a_j e_j$ are those with $a_0 =0$. So, there is no non-trivial, closed, invariant subspace which, as in \eqref{eq11}. has a distance between $0.3$ and $0.7$ to $e_0$. But, for such weighted shifts, it may require very large $\|\ell'\|_2$ in order to make 
$$\ell'(T) (e_0 + 10 e_n) - e_0 < 0.7.$$
So we may produce many ``$(\varepsilon \theta)_n$-almost cyclic vectors'' but we will not get norm convergence to a non-cyclic vector.

If, throughout our MC process, there exists a sequence of positive real numbers $(\delta_k)$ such that 
$$\left| \langle T^j y'_n, y'_n \rangle \right| > \delta_k \|y'_n\|^2$$
for some $j \ge k $ (which may depend on $n$) then our MC will have $y'_n$ converging to a non-cyclic vector. If there is no such sequence $(\delta_k)$ we will combine our MC with other arguments. We either get a vector $y'$, $\|y'\| =1$, and $m \ge 1$ such that  
$$\langle T^j y', y' \rangle = 0 \text{ for every } j \ge m,$$
i.e., $T^my'$ is non-cyclic, or a vector $y'$, $\|y'\|=1$, and an $m \ge 1$ such that 
$$\langle T^jy', x_0 \rangle = 0 \text{ for } j \ge m,$$
i.e., 
\begin{equation}\label{eq12}
T^m y' \quad \text{is non-cyclic.}
\end{equation}

In our MC we will construct a sequence $(y_n)$ by the formula
\begin{equation}\label{eq13}
y_{n+1} = y_n + \sum_{j \ge 0} r_j T^j y_n.
\end{equation}
By using \eqref{eq13} instead of the mere general $y_{n+1} = y_n + z$ we will be able to use elementary Fourier Analysis in order to give the sequence $(y_n)$ the desired properties.

From equation \eqref{eq05} we get
$$\ell'(T) = C \left( I + C V_{y'}^{*} V_{y'}\right)^{-1} V_{y'}^{*} x_0 \quad \text{ with } C = 1/ C'.$$
Observing that 
$$\left( I + C V_{y'}^{*} V_{y'}\right)^{-1} V_{y'}^{*} = V_{y'}^{*} \left( I + C V_{y'} V_{y'}^{*} \right)^{-1} $$
we get
$$V_{y'}\ell' = V_{y'} V_{y'}^{*} C \left( I + C V_{y'} V_{y'}^{*}\right)^{-1} x_0.$$
Next, putting $y = C^{\frac{1}{2}} y'$ and $\ell = C^{-\frac{1}{2}} \ell'$ we get
\begin{equation}\label{eq14}
V_{y'}\ell' = V_y \ell = V_y V_y^{*}  \left( I + V_{y}V_{y}^{*} \right)^{-1} x_0
\end{equation}
So, if follows that $V_y V_y^{*}  \left( I + V_{y} V_{y}^{*} \right)^{-1} x_0$ is the element that we get when we move $y$ within the distance $\|(I + V_y V_y^{*})^{-1} x_0\|$ to $x_0$, with an $\ell(T)$ with $\ell$ of minimal $\|\cdot \|_2$-norm. From \eqref{eq14} we see that, when moving $y$ within some distance from $x_0$, with $\ell$ of minimal $\|\cdot \|_2$-norm, we have 2 representations: one is $\ell'(T) y'$ with $\|y'\| \approx 1$ and the other is $V_y V_y^{*}  \left( I + V_{y} V_{y}^{*}\right)^{-1} x_0$ with $y = C^{\frac{1}{2}} y'$.

We see that $V_{y} V_{y}^{*}$ is a compact, self-adjoint operator for which $\langle x, V_{y} V_{y}^{*} x \rangle \ge 0 \, \forall x.$  We will now, for short, put 
\begin{equation}\label{eq15}
\left(I+V_{y} V_{y}^{*}\right)^{-1} = \VV, \text{ so} \left(I+V_{y} V_{y}^{*}\right)^{-1}  = \VV x_0 \text{ and } x_0 - \ell(T)y = \VV x_0.
\end{equation}

%%%%% INSERT 1
%%%%% %%%%% 
When moving $y$ within distance $\varepsilon$, $0<\varepsilon <1$, from $x_0$ by an $\ell (T) =  \sum_{j \ge 0} a_j T^j$ with $\sum_{j \ge 0} |a_j|^2$ minimal, we arrive at the same point $\ell (T) y = x_0 - \VV x_0$, whether we move $y$ or  $y \cdot e^{ir}$, $r \in \mathbb{R}$, or whether we use the operator $T$ or the operator $ e^{is} T $, $s \in \mathbb{R}$.

This affects $a_m = \langle \VV x_0, T^m y \rangle$ and $b_m = \langle [ \,\, ]^{-2} x_0, T^m y \rangle$ and below without loss of generality we can assume the values of the integrals in \eqref{eq28primeprime}--\eqref{eq30primeprime} below to be real numbers even if we just consider the parts of the integrals with $\displaystyle \sum_{j \ge 0} r_j e^{i j \theta}$.
%%%%% %%%%% 
%%%%%% %%%%% 

We get, by \eqref{eq04}, 
\begin{equation}\label{eq16}
\|\ell\|^2_2 = \sum_{j \ge 0} \left| \langle \VV x_0, T^j y\rangle\right|^2 = \langle V_y^* \VV x_0, V_y^* \VV x_0 \rangle = \langle V_y V_y^* \VV x_0, \VV x_0\rangle = \varepsilon \theta > 0.
\end{equation}
So, with $\ell(T) y = \sum_{j \ge 0} a_j T^jy$ we get 
\begin{equation}\label{eq17}
a_j = \langle \VV x_0, T^j y\rangle \text{ and } \sum_{j \ge 0} |a_j|^2 = \varepsilon \theta .
\end{equation}
Moreover, with $\ell_{\varepsilon}(T)$ having minimal $\|\cdot\|_2$-norm for $\ell$ is with $\|x_0 - \ell(T) y\| \le \varepsilon$, we get from \eqref{eq14} that  
\begin{equation}\label{eq18}
\varepsilon \rightarrow \ell_{\varepsilon} (T) \text{ is a continuous function.}
\end{equation}

We will now divide the operators on $H$ that we consider ($T$ is one to one, $\mathcal{R}(T)  \neq H$, $\overline{\mathcal{R}(T)}=H$, $0 \in \sigma (T)$,  $\|T\|_{\text{op}} = 10^{-20}$) into two disjoint classes that we shall refer to as Type 1 and Type 2.
\bigskip 

\noindent \textbf{\underline{Type 1.}}  We say that $T$ is of type 1 if there is $u_0 \in H, $ $\|u_0\| = 1$ such that for every integer $n\ge 1$, there is $\delta_n > 0$ such that for every $y$ with $\displaystyle \langle \frac{y}{\|y\|}, u_0 \rangle \ge \frac{1}{100}$ there is $j \ge n$ with
\begin{equation}\label{eq19}
	\left| \langle T^j y, y\rangle\right| \ge \delta_n \|y\|^2.
\end{equation}

\bigskip 

\noindent \textbf{\underline{Type 2.}}  We say that $T$ is of type 2 if  for every  $u_0 \in H, $ $\|u_0\| = 1$, there is an integer $m \ge 1$ such that for every $\delta > 0$ there is $y$ with 
$\displaystyle \langle \frac{y}{\|y\|}, u_0 \rangle \ge \frac{1}{100}$ such that
$$	\left| \langle T^j y, y\rangle\right| \le \delta \|y\|^2$$
for every $j \ge m$.

For every operator of type 1 we will construct by (MC) an invariant subspace of the type described in \eqref{eq11}. For operators of type 2 we proceed as follows:  Choose $u_0 \in H$ and $m$ as in the definition of type 2 operators. Let $(\delta_n)$ be a sequence of positive numbers, with $\delta_n \rightarrow 0$. For every $\delta_n$, choose $y_n$ with $\displaystyle \langle \frac{y}{\|y\|}, u_0 \rangle \ge \frac{1}{100}$ such that 
$$	\left| \langle T^j y_n, y_n \rangle\right| \le \delta_n \|y_n\|^2.$$
for every $j \ge m$.

Let a subsequence of $y_n / \|y_n\|$ converge, weakly or in norm, to $y'_{\infty}$ with 
$$\displaystyle \langle y'_{\infty},u_0\rangle \ge \frac{1}{100}.$$ 

If the convergence is in norm then $\displaystyle \langle T^j y'_{\infty}, y'_{\infty} \rangle = 0$ for all $j \ge m$, so $T^m y'_{\infty}$ is non-cyclic. If the convergence is just weak, we put 
$$\frac{y_n}{\|y_n\|} = \alpha_n y'_{\infty} + s_n$$
with $|\alpha_n| \ge 1/100$. Then

\begin{equation}\label{eq20}
\langle T^j \alpha_n y'_{\infty},\alpha_n y'_{\infty} \rangle + \langle T^j s_n, s_n\rangle \longrightarrow 0 \text{ for every } j \ge m
\end{equation}

Below we will use this to construct an invariant subspace.

Equation \eqref{eq20} gives, for fixed $p$ and $q$, 

\begin{equation}\label{eq21}
\langle p(T) \alpha_n y'_{\infty}, (T^*)^m q(T^*) \alpha_y y'_{\infty} \rangle + \langle p(T)s_n, (T^*)^m q(T^*) s_n \rangle \longrightarrow 0.
\end{equation}

For operators of type 1 we will choose $x_0$ to be 
$$x_0 = \frac{\sqrt{3}}{2} u_0 + \frac{1}{2} u_1,$$
where $u_0$ is from definition of type 1 and $u_1$ is chosen such that 
$$\| T^* u_1\| < (\varepsilon \theta)_0,$$
with $(\varepsilon \theta)_0$ chosen below, $\langle u_1, u_0 \rangle = 0$, $\|u_1\| = 1$, $(\varepsilon \theta)_0$ is the $(\varepsilon \theta)$ from which we below start our (MC).

For operators of type 2, both $u_0$ and $u_1$, $\|u_0\| = \|u_1\| = 1$, $\langle u_1, u_0\rangle = 0$, can be chosen arbitrarily.

In many estimates below we can keep in mind that the $(\varepsilon \theta)$'s and the functions thereof, like $\gamma (\varepsilon \theta)$, will be of a smaller order of magnitude than the $(\delta_k)$'s in the definition of type 1.

In our (MC) in $\ell (T) y_n = \sum_{j \ge 0} a_j T^j y_n$ it will be important to have $(a_0)$ dominating,
\begin{equation}\label{eq22}
	|a_0|^2 \ge 100 \sum_{j \ge 1} |a_j|^2
\end{equation}
and we will now make preparations to arrange this. Lemma 2 shows how we can get, from a general
$$\ell'(T) y' = \sum_{j \ge 0} a_j' T^j y',$$
a
$$\ell^{\prime \prime} (T) = \sum_{j \ge 0} a''_j T^j$$
with
$$\ell'(T)y' = y_1', \quad \ell^{\prime \prime} (T) y_1' = \sum_{j \ge 0} a''_j T^j y_1' \, \text{ and } \, |a''_0| \, \text{ dominating in } \, \sum_{j \ge 0} a''_j T^j$$
such that we can start our (MC) from 
$$ \ell^{\prime \prime} (T) y_1' = \sum_{j \ge 0} a''_j T^j y_1' .$$

It will also be useful to have estimates from below on $\sum_{j \ge 1}|a_j''|^2$ compared to $|a''_0|^2$. Except for a singular case, when we have an invariant subspace anyhow we will arrange that: when we start (MC) we have estimates both from below and above of $\displaystyle \frac{\sum_{j \ge 1}|a_j|^2}{|a_0|^2}$ in $\ell (T) = \sum_{j \ge 0} a_j T^j$.

With $\displaystyle \sum_{j \ge 0} |a_j|^2 = \varepsilon \theta$ from \eqref{eq16} we get, from \eqref{eq22}, $$\displaystyle (\varepsilon \theta)^{1/2} > |a_0| > 0.1 (\varepsilon \theta)^{1/2}$$ and, since, 
$$a_0 y + \sum_{j \ge 1} a_j T^j y = \ell (T) y$$
we get
\begin{equation}\label{eq23}
\frac{1}{20} (\varepsilon \theta)^{-1/2} < \| y \| < (\varepsilon \theta)^{-1/2}.
\end{equation}

For the next Lemma 1 we put, as above, $x_0 = \sqrt{3}/2 u_0 + u_1$. We assume $T$ to be of type 1. We consider, for $0.3 < \varepsilon \le 0.5,$

$$\| \ell'_{\varepsilon} \|_2 = \text{inf } \left\{ \|\ell'\|_2 , \|x_0 - \ell'(T y_0' \|)\right\}  \le \varepsilon, \quad y_0' = \frac{\sqrt{3}}{2} u_0.$$
\bigskip

\noindent \textbf{Lemma 1.} For some $\varepsilon > 0$, 
$$0.5 \ge \varepsilon > 0.5 - 10^{-5} (\varepsilon \theta)_0$$
we have, for the corresponding $\ell'_{\varepsilon}(T)$ and $\varepsilon \theta,$
$$(\varepsilon \theta) \le \frac{1}{2} \, 10^{-5} (\varepsilon \theta)_0.$$

\noindent {\sc Proof of Lemma 1}.\\
\noindent We have $\| x_0 - I y_0'\| = 0.5$, so $\ell'_{0.5} \le 1$. On the other hand, since $\|T\|_{\text{op}} = 10^{-20}$, $\|\ell'_{0.5}\|_2 > 0.8$. If 
$$\|\sum_{j \ge 0} a_j T^j y'_0 - x_0 \| \le 0.5 - 10^{-5} (\varepsilon\theta)_0,$$
then
$$\langle \sum_{j \ge 1} a_j T^j y'_0 , u_1 \rangle \ge 10^{-5} (\varepsilon\theta)_0,$$
and so 
\begin{equation}\label{eq24}
\displaystyle \text{max}_{j \ge 1} |a_j| \ge 10^{20-5-1} = 10^{14} > \|\ell'_{0.5}\| \cdot 10^{13}. 
\end{equation}

If we assume that the $(\varepsilon \theta)$ from $\ell'_{\varepsilon}$, say $(\varepsilon \theta)'_{\varepsilon}$ is $ > \frac{1}{2} 10^{-5} (\varepsilon \theta)_0$ for every $\varepsilon$, $0.5 \ge \varepsilon > 0.5 - 10^{-5} (\varepsilon \theta)_0$, we get from \eqref{eq10} that
$$\left\|   \ell'_{\varepsilon - \frac{\delta}{10} \cdot 10^{-5} (\varepsilon \theta)_0}\right\|_2 < \|\ell'_{\varepsilon}\|_2 \, (1 + \delta)$$
for $\delta$ sufficiently small. Applying this $10/\delta$ times, starting with $\varepsilon = 0.5$ and finishing with $\varepsilon = 0.5 - 10^{-5} (\varepsilon \theta)$ we get
$$\left\| \ell'_{0.5 - 10^{-5} (\varepsilon \theta)_0}\right\| < (1+\delta)^{10/\delta} \|\ell'_{0.5}\|_2 < 10^{10} \|\ell'_{0.5}\|_2$$
which contradicts \eqref{eq24}. This proves Lemma 1. \hfill $\Box$

We shall now provide a series of new lemmas with their corresponding proofs throughout the forthcoming pages.

In the next lemma (Lemma 2), we show that from a situation with $\ell (T) = \sum_{j \ge 0} a_j T^j$ and some small $(\varepsilon \theta)$ we get to a situation with a dominating $|a_0|$ in $\sum_{j \ge 0} a_j T^j$. This is initially at the cost of possibly getting a larger $\varepsilon \theta$. But, since, form there, we can use MC, we can get the $(\varepsilon \theta)_n$ to tend to $0$.
\bigskip

\noindent \textbf{Lemma 2.} Put $\ell (T) y ) y'_1$, $\langle y_1', x_0 - y_1' \rangle = \varepsilon\theta$, and $\|\ell\|_2$ minimal. Choose $s$ such that $s y_1' = y_1$, we have
$$\|x_0 - \sum_{j \ge 0} a_j T^j y_1 \| = \|V_{y_1} V^{*}_{y_1} \VV x_0 - x_0 \| = \|x_0 - y_1'\| .$$
There is a function $\gamma = \gamma (\varepsilon \theta) = \gamma (\varepsilon \theta, u_0, T)$ defined for $\varepsilon\theta \le 10^{-4}$ such that, for $\| \sum_{j \ge 0} a_j T^j y_1' - x_0\| =  \|x_0 - y_1'\| ,$ $$a_0^2 \ge \left[ 1 - \gamma (\varepsilon \theta)\right], \quad \gamma(\varepsilon \theta) \rightarrow 0 \text{ as } (\varepsilon \theta) \rightarrow 0.$$
\bigskip

\noindent \textbf{Lemma 3.} Under the same assumptions as in lemma 2, with $\langle y_1', x_0 - y_1'\rangle = \varepsilon \theta$ we have
$$\langle \sum_{j \ge 0} a_j T^j y_1, x_0 - \sum_{j \ge 0} a_j T^j y_1\rangle \le 3 \left[ \gamma (\varepsilon\theta)\right]^{1/2}.$$

We also have:
\bigskip

\noindent \textbf{Lemma 4.}  Under the same assumptions as in lemma 2,
$$\|y_1' - \sum_{j \ge 0} a_j T^j y_1\|\le \left[ \gamma (\varepsilon\theta)\right]^{1/2}.$$ 

Lemma 4 provides the following:
\bigskip

\noindent \textbf{Lemma 5.} 
$$\langle \sum_{j \ge 0} a_j T^j y_1, x_0 - \sum_{j \ge 0} a_j T^j y_1\rangle \le 3 \left[ \gamma (\varepsilon\theta)\right]^{1/2}.$$

\bigskip
\noindent {\sc Proof of Lemma 2}.\\
Let $K = 10^{20}$. We have $\displaystyle \left|\langle T^m y'_1, x_0 - y_1'\rangle \right| \le \varepsilon\theta$ for all $m \ge 0$ and for 
$$m = \left[ 1+ \log_{K} \left(\frac{1}{\varepsilon\theta}\right)\right]$$
we get 
$$\left\| T^{m+r} y_1'\right\| \le (\varepsilon\theta) \cdot \frac{1}{K^r}.$$
This gives
$$\left| \sum_{j \ge 0} c_j T^j y_1', x_0 - y_1'\right| \le \left( 2 \log_{K} \left(\frac{1}{\varepsilon\theta}\right) \right) \cdot \varepsilon\theta$$
if $|c_j| \le 1$ for every $j$.

With 
$$\sum_{j \ge 0} a_j T^j y_1 = \sum_{j\ge0} c_{1j} T^j y_1' = (1-d) y_1' + \sum_{j \ge 1} c_{1 j} T^j y_1',$$
the minimality of
$$|(1-d)^2| + \sum_{j \ge 1} |c_{1 j}|^2$$
for the distance $\|x_0 - y_1'\|$ to $x_0$ gives
$$|c_{1 j}| \le 2 |d|^{1/2}$$
for every  $j \ge 1.$ 

Now, put 
$$h (|d|) = \text{inf } \left\{ \|d y_1' - \sum_{j \ge 1} c_j T^j y_1' \|, |c_j|\le 1\right\}$$
for every $j \ge 1$.

We first observe that $h (|d|) > 0$. Otherwise, we could find a sequence $y'_{1 n}$, $\langle y'_{1 n}, u_0\rangle > 1/100$, weakly convergent, such that for $w-\text{lim} y'_{1 n} = y'_{1 \infty}$ we have 
$$d y'_{1 \infty} = \sum_{j \ge 1} c_j T^j y'_{1 \infty}$$
which contradicts that $y'_{1 \infty}$ is cyclic and $T$ is not invertible. Moreover, for $\kappa$  being small and positive we have
\begin{align*}
(1+\kappa) h(|d|)&= \text{ inf } \left\{\|d (1+ \kappa ) y'_1 - \sum_{j \ge 1}c_j T^j y'_1\| \, , \, |c_j| \le 1 + \kappa \right\}\\
& < \text{ inf } \left\{ \|d (1+\kappa) y'_1 - \sum_{j \ge 1}c_j T^j y'_1\| \, , \, |c_j| \le 1 \right\}\\
& = h ((1+\kappa)|d|) 
\end{align*}
so $$\frac{h(|d|)}{|d|}$$ is increasing in $|d|$.

We get
\begin{align*}
	\|x_0 - y_1'\|^2 & = \left\| x_0 - \left(  (1-d)y_1' + \sum_{j\ge 1} c_{1 j} T^j y_1'  \right) \right\|^2 \\
	& = \|x_0 - y_1'\|^2 + \|d y_1' - \sum_{j\ge 1} c_{1 j} T^j y_1' \|^2\\
	& - 2 |\langle d y_1' - \sum_{j \ge 1} c_{1 j} T^j y_1', x_0  - y_1'\rangle | \\
	& \ge \|x_0 - y_1'\|^2 +h^2(|d|) - \left( 4 \log_K \left(\frac{1}{\varepsilon\theta}\right)\right) \cdot \varepsilon\theta .
\end{align*}
This gives
$$h(|d|) \le 2 \left[\log_K \left(\frac{1}{\varepsilon\theta}\right)\right]^{1/2} $$
and, with $h(|d|)$ strictly increasing, we get
$$|d| \le h^{-1} \left(  2 \left[\log_K \left(\frac{1}{\varepsilon\theta}\right)\right]^{1/2}  \right) = \gamma (\varepsilon \theta).$$

Before continuing after $\ell_n(T)y_n = y_1'$ and $\displaystyle \sum_{j \ge 0} a_j T^j y_1'$ we shall arrange to have an estimate from below of $\sum_{j \ge 1} |a_j|^2$ in $\sum_{j \ge 0} a_j T^j y_1'$. We distinguish between two cases. We first put $\|T\|_{\text{op}} = \frac{1}{K}$, with $K = 10^{20}$.

\bigskip

\noindent \textbf{\underline{Case I.}} $\displaystyle |\langle T^j y_1', x_0 - y_1'\rangle| \ge (\VE)^4$ for some $j$.\\ 
\noindent Then  $$\|x_0 - \left[ (1- (\VE)'') y_1' + (\VE)^6 \, T^j y_1' \right]\| \le \|x_0 - y_1'\| + (\VE)'' - (\VE)^{10} + (\VE)^{12} < \|x_0 - y_1'\|$$
and with
$$\left( 1- (\VE)''\right)^2 + \left( (\VE)^6 \right)^2 < 1$$
we see that, for $\text{inf } \sum_{j \ge 0} |a_j|^2$ in $\displaystyle \| \sum_{j \ge 0} a_j T^j y_1' - x_0\| \le  \|x_0 - y_1'\|$ we have 
\begin{equation}\label{eq25}
	 \sum_{j \ge 1} |a_j|^2 \ge (\VE)^{15}.
\end{equation}
\bigskip

\noindent \textbf{\underline{Case II.}} $\displaystyle |\langle T^j y_1', x_0 - y_1'\rangle| < (\VE)^4$ for all $j$.\\  
\noindent We have 
$$\sum_{j \ge 1} \left| \langle c_j T^j y_1', x_0 - y_1'\rangle\right| < 2 |\log_{K} (\VE)| \cdot (\VE)^4$$
if $|c_j| \le 1$ for all $j$.\\

For $\displaystyle \delta = \frac{\VE}{10}\,$ let $\|\ell' (T)\|_2$ be minimal for 
\begin{equation}\label{eq26}
	\|x_0 - \ell' (T) y_1'\| \le \|x_0 - (1 + \delta) y_1'\|.
\end{equation}

We get 
\begin{equation}\label{eq27}
\|\ell'(T) y_1' - (1+\delta) y_1'	\| < (\VE)^2.
\end{equation}

To see \eqref{eq27} let us observe that in
$$\ell'(T) y_1' = \sum_{j \ge 0} a_j T^j y_1'$$
we cannot have $|a_0| \ge (1+ \delta)$. On the other hand, if 
$$|a_0| < 1 + \delta - (\VE)^2$$
then, to fulfill \eqref{eq24} we would need $|a_j| \ge 1$ for some $j \ge 1$.

Equation \eqref{eq27} gives 
$$\langle \ell'(T) y_1', x_0 - \ell' (T) y_1' \rangle < \VE \left( 1 - \frac{1}{20}\right),$$
so we can now use Lemma 2 again, letting $\ell'(T)y_1'$ of \eqref{eq27} playing the role of $\ell (T) y$ in Lemma 2, and letting
$$\ell'(T) y_1' = y_1'' ,$$ 
for some $y_1''$, play the role of $y_1'$ in Lemma 2 and letting 
$$(\VE)' = \langle y_1'', x_0 - y_1''\rangle = \langle \ell'(T) y_1', x_0 - \ell'(T) y_1'\rangle < \VE \left( 1 - \frac{1}{20}\right)$$
play the role of $\VE$ in Lemma 2.

Thus, every time we apply Lemma 2 we either get Case I and we are done, or we get Case II and pass to a smaller $\VE$, diminished by a factor $(1 - \frac{1}{20})$, and not moving more than $(\VE)^2$ by \eqref{eq27}. So, if Case II happens every time we obtain convergence to a non-cyclic vector.

%%%%%%%% INSERT 2
%%%%%%%%%%%%%%
Now we will study the change of $\VV x_0$ when we make the change $y \rightarrow y (1+\delta)$ and when we make the change  $ y \rightarrow y + \alpha T^j y,$ $\alpha$ a small complex number.

If $A$ is an invertible operator the first order change of $A^{-1}$ when we make the change  $A \rightarrow A+B$ is $A^{-1} \rightarrow A^{-1} - A^{-1} B A^{-1}$, $y \rightarrow y (1+\delta)$  gives in first order
$$\VV x_0 \rightarrow \VV x_0 - 2 \delta \VV V_y V_y^* \VV x_0 = \VV x_0 - 2 \delta ([\,\,]^{-1} - [\,\,]^{-2}) x_0.$$
With $\langle \VV x_0, T^m y \rangle = a_m$ and $\langle [\,\,]^{-2}x_0, T^m y \rangle = b_m$ we get  
\begin{align*}
\sum a_m T^m y & = x_0 - \VV x_0\\
2 \delta \sum b_m T^m y &= 2 \delta (\VV - [\,\,]^{-2}) x_0
\end{align*}
%%%%%%%%%%%%%%
%%%%%%%%%%%%%%
\textbf{Remark.} We can think of Case II happening  every time as a singular situation. Every time we move $y$ closer to $x_0$ by a new minimal $\ell (T)$, say $\ell_{1 \delta}(T),$ $\ell_{1 \delta}(T)$ will just be near to $\ell(T)  + 2 \delta I$ and $\ell_{1 \delta}(T) y$ will be near to $\ell(T) y + 2 \delta y$.  However, $\ell_{1 \delta}(T)y$ will also be near to 
$$x_0 - \VV x_0 + 2 \delta (\VV - [\, \,]^{-2})x_0 = \sum a_m T^m y + 2 \delta \sum b_m T^m y.$$
So this shows that $b_m$ is near $0$ for $m \ge 1$ and $(\VV - [\, \,]^{-2})x_0$ is almost collinear with $y$. Equation \eqref{eq25} can be sharpened by a similar argument as above with Case I and Case II, with the conclusion that we either have a non-cyclic vector or there are at least 4 different $j$'s, bigger than $1$, say $j_1, j_2, j_3, j_4$ such that $|a_{j_1}|$, $|a_{j_2}|$, $|b_{j_3}|$, $|b_{j_4}|$ are all bigger than $(\VE)^{30}$.
%%%%%%%%%%%%%%
%%%%%%%%%%%%%%
%Now we will first study the first order change of $\VV x_0$ in $\alpha$, $\alpha$ being a small complex number when we change $y$ to $y + \alpha T^j y$.
%
%If $A$ is an invertible operator, the first order change of $A^{-1}$ when we make the change $A \mapsto A+B$ is $A^{-1} \mapsto A^{-1} - A^{-1} B A^{-1}$, 
\bigskip \bigskip

The first order in $\alpha$, $ y \mapsto y + \alpha T^j y$ gives 
$$V_y V_y^* \mapsto V_y V_y^* + \alpha V_{T^jy} V_y^* + \overline{\alpha} \, V_y V^*_{T^jy},$$
$$\VV x_0 \mapsto \VV x_0 - \VV \, \left( \alpha V_{T^jy} V_y^* + \overline{\alpha} \, V_y V^*_{T^jy}\right)\, \VV x_0.$$

The first order change of $\langle \VV x_0, x_0\rangle$ is given by
$$\langle \VV x_0, x_0\rangle \mapsto \langle \VV x_0, x_0\rangle - \langle \VV \, \left( \alpha \, V_{T^jy} V_y^* + \overline{\alpha} \, V_y V^*_{T^j y}\right) \VV \,x_0, x_0 \rangle .$$
With $\displaystyle \langle \VV x_0, T^m y\rangle = a_m$ the change becomes
\begin{equation}\label{eq28}
- \sum_{m \ge 0} \left( \alpha a_m \overline{a_{m+j}} + \overline{\alpha} a_{m+j} \overline{a_m}\right).	
\end{equation}

With $b_m = \langle [\, \,]^{-2} x_0, T^m y\rangle$ the change becomes

\begin{equation}\label{eq29}
	- \sum_{m \ge 0} \left( \alpha a_m \overline{b_{m+j}} + \overline{\alpha} a_{m+j} \overline{b_m}\right).	
\end{equation}

The $a_0$ coefficient in $\displaystyle \sum_{m\ge 0} a_m T^m y$ is $\langle \VV x_0, y\rangle$, which will in first order in $\alpha$ change to
$$\langle \VV x_0 - \VV \left( \alpha \, V_{T_y^j} \,V_{y^*} + \overline{\alpha} V_y \, V^*_{T^jy} \right) \VV x_0, y + \alpha \, T^j y \rangle.$$

With $\kappa_m = \langle \VV y, T^m y\rangle$ this becomes
\begin{equation}\label{eq30}
	\langle \VV x_0, y \rangle + \overline{\alpha} \, a_j \, (1 - \kappa_0) - \left( \sum_{m \ge 0} \alpha \, a_m \overline{\kappa_{m+j}} + \overline{\alpha} \, \sum_{m \ge 1} a_{m+j} \overline{\kappa_m}\right).
\end{equation}

Put 
$$f(e^{i \theta}) = \sum_{m \ge 0} a_m e^{i m \theta}, \quad g(e^{i \theta}) = \sum_{m \ge 0} e^{i m \theta}, \quad \kappa(e^{i \theta}) = -1 + \sum_{m \ge 0} \kappa_m \, e^{i m \theta}.$$
We see that the changes from \eqref{eq28}, \eqref{eq29}, \eqref{eq30} are

\begin{equation}\label{eq28prime}\tag{28$^\prime$}
-\frac{1}{2\pi} \int_{0}^{2\pi} \left(\alpha \, e^{i j \theta} + \overline{\alpha} \, e^{- i j \theta}\right)\, f(e^{i\theta}) \, \overline{f(e^{i\theta})} \, d\theta
\end{equation}

\begin{equation}\label{eq29prime}\tag{29$^\prime$}
-\frac{1}{2\pi} \int_{0}^{2\pi} \left(\alpha \, e^{i j \theta} + \overline{\alpha} \, e^{- i j \theta}\right)\, f(e^{i\theta}) \, \overline{g(e^{i\theta})} \, d\theta
\end{equation}

\begin{equation}\label{eq30prime}\tag{30$^\prime$}
	-\frac{1}{2\pi} \int_{0}^{2\pi} \left(\alpha \, e^{i j \theta} + \overline{\alpha} \, e^{- i j \theta}\right)\, f(e^{i\theta}) \, \overline{\kappa (e^{i\theta})} \, d\theta
\end{equation}

The changes given by 
$$y \mapsto y + \sum_{j \ge 0} r_j T^j y$$ will be in first order

\begin{equation}\label{eq28primeprime}\tag{28$^{\prime\prime}$}
	-\frac{1}{2\pi} \int_{0}^{2\pi} \left[ \sum_{j \ge 0} \left(
	r_j e^{i j \theta} + \overline{r_j} e^{- i j \theta} \right) \right]  \, f(e^{i\theta}) \, \overline{f(e^{i\theta})} \, d\theta
\end{equation}

\begin{equation}\label{eq29primeprime}\tag{29$^{\prime\prime}$}
	-\frac{1}{2\pi} \int_{0}^{2\pi} \left[ \sum_{j \ge 0} \left(
	r_j e^{i j \theta} + \overline{r_j} e^{- i j \theta} \right) \right]  \, f(e^{i\theta}) \, \overline{g(e^{i\theta})} \, d\theta
\end{equation}

\begin{equation}\label{eq30primeprime}\tag{30$^{\prime\prime}$}
	-\frac{1}{2\pi} \int_{0}^{2\pi} \left[ \sum_{j \ge 0} \left(
	r_j e^{i j \theta} + \overline{r_j} e^{- i j \theta} \right) \right]  \, f(e^{i\theta}) \, \overline{\kappa (e^{i\theta})} \, d\theta
\end{equation}

When we take a step

$$y \mapsto y + \sum_{j \ge 0} r_j T^j y$$

in MC, the change in $\VV x_0$ will be denoted $ch\VV x_0$.

We want $\VE$ to decrease which means

\begin{equation}\label{eq31}
\langle ch\VV x_0, x_0 - 2 \VV x_0\rangle < 0 \text{ in first order.}
\end{equation}

We want 
$$\left| \langle ch\VV x_0, x_0\rangle  \right| \, \text{ and } \, \langle ch\VV x_0, \VV x_0\rangle $$ not to be bigger than
$$10 |\langle ch\VV x_0, x_0 - 2 \VV x_0\rangle |$$
in order to ensure that 

\begin{equation}\label{eq32}
	\|\VV x_0\| \, \text{ does not tend to } \, 0 \text{ or }\, 1.
\end{equation} 

And we also want $|a_0|^2$ to follow the size of  $(\VE)$ such that \eqref{eq25} and Lemmas 3, 4, and 5
\begin{equation}\label{eq33}
	1 - \frac{1}{2} (\VE)^{15} > \frac{|a_0|^2}{\VE} > 1 - \left[ (3 \gamma^{1/2}) \circ (3 \gamma^{1/2} \right]\, \VE
\end{equation}
We say that a finite sequence $Z_n$, $n_0 \le n \le n_1$ of $m$ elements in $H$,
$$Z_n = \left( z_{1 n}, z_{2 n}, \ldots, z_{m n} \right)$$
is $\VE$-linearly independent in an interval $t \le \VE \le u$ if there is a positive continuous function $s(\VE)>0$ for $\VE > 0$ such that
$$\left\| \sum_{i=1}^{m} a_{i n} \, z_{i n}\right\| \ge s(\VE) \, \text{max}_i \, |a_{i n}|$$
for  every $n$.

If the sequence of the three functions  $f(e^{i \theta})$, $g(e^{i \theta})$, and $\kappa (e^{i \theta})$ created from the sequence of $y_n$'s, $n_0 \le n \le n_1$, is $\VE$-linearly independent in the interval of $\VE$'s created by the $y_n$'s, then we can fulfill \eqref{eq31}--\eqref{eq33}. To see this, we observe that  

\begin{equation}\label{eq34}
	\begin{aligned}
&\text{the size of the integrals in \eqref{eq28primeprime}--\eqref{eq30primeprime} are bounded from} \\
& \text{ below by } (\sum_j |r_j|) \times (\text{some fixed positive function of  } \VE ).
\end{aligned}
\end{equation}

In order to verify \eqref{eq34}, see, e.g. \cite{BeauzamyEnflo} and observe that the size of the $j$'s are bounded by some fixed function of $(\VE)$.

To have $f(e^{i \theta})$,  $g(e^{i \theta})$, $\kappa (e^{i \theta})$ $\VE$-linearly independent for a finite sequence of $y_n$'s and $(\VE)$ in some given interval it is sufficient to have corresponding $\VE$-linear independence for 
\eqref{eq35}--\eqref{eq37} below:

\begin{equation}\label{eq35}
	V_y V_y^* \VV x_0 = \sum_{j \ge 0} \langle \VV x_0, T^j y \rangle T^jy = \sum_{j\ge 0} a_j T^j y
\end{equation}

\begin{equation}\label{eq36}
	V_y V_y^* [\, \,]^{-2} x_0 = \sum_{j \ge 0} \langle  [\, \,]^{-2} x_0, T^j y \rangle T^jy = \sum_{j\ge 0} b_j T^j y
\end{equation}
 and
 \begin{equation}\label{eq37}
\begin{aligned}
 	-y + \sum_{j \ge 0}V_y V_y^* [\, \,]^{-1} y &= (-1+\kappa_0) y + \sum_{j \ge 1} \langle  [\, \,]^{-1} y, T^j y \rangle T^jy \\
 	& = (-1+\kappa_0) y + \sum_{j\ge 1} \kappa_j T^j y = \VV y.
\end{aligned}
 \end{equation}

 So, if \eqref{eq35}--\eqref{eq37} are fulfilled we can use MC and get down to arbitrarily small $\VE$'s and keep
 $$0.3 < \langle \VV x_0, x_0 \rangle < 0.7.$$

However, this is not enough to produce non-cyclic vectors, and we will get back to that below. First we will show that even if we do not have the $\VE$-linear independence in \eqref{eq35}--\eqref{eq37} we can use MC to get down to arbitrarily small $\VE$'s. We recall that we are studying operators of type 1. 

If $\VE$-linear independence is not fulfilled in \eqref{eq35}--\eqref{eq37} we can solve the inequality

\begin{equation}\label{eq38}
	\left\| A \,  V_y \, V_y^* \,  \VV \,  x_0 + B  \,  V_y \,  V_y^* \,  [\, \,]^{-2} \,  x_0 - \VV y\right\| \le (\VE)^{100}.
\end{equation}

We can apply $(I + V_y V_y^*)$, which has operator norm less than $ \frac{1}{(\VE)}$ for each term and get
$$\|A V_y V_y^* x_0 + B V_y V_y^* \VV x_0 - y\| < (\VE)^{98}$$
 which gives 
 
 \begin{equation}\label{eq39}
 \|A V_y V_y^* (x_0 - \VV x_0) + (A+B) V_y V_y^* \VV x_0 - y\| < (\VE)^{98} .
 \end{equation}
 
To solve \eqref{eq39} we get from Lemmas 2--5 that

\begin{equation}\label{eq40}
	\begin{aligned} 
&\left\| \frac{V_y V_y^* \VV x_0}{\|V_y V_y^* \VV x_0\|} - \frac{y}{\|y\|}\right\| \le \gamma' (\VE)  \\
&\text{ where } \gamma' = (3 \gamma^{1/2}) \circ (3 \gamma^{1/2}) \text{ from \eqref{eq33} and } \gamma'(\VE)  \stackrel{\VE \rightarrow 0}{\longrightarrow} 0.
	\end{aligned}
\end{equation}
We have
\begin{equation}\label{eq41}
\left\| \frac{V_y V_y^* (x_0 - [\, \,]\,  x_0)}{\|V_y V_y^* (x_0 - [\, \,]\,  x_0)\|} - \frac{y}{\|y\|}\right\| > \frac{1}{100}\, \delta_1  \|T^j y/\|y\|\| \ge \frac{1}{100} \delta_1^2,  \text{ since } j < |\log_K \delta_1|
\end{equation}
 $T$ is of type 1 (see \eqref{eq19}).
 
 The same estimate also works for
 \begin{equation}\label{eq42}
 	\left\| \frac{V_y V_y^* x_0}{\|V_y V_y^* x_0\|} - \frac{y}{\|y\|}\right\| > \frac{1}{100}\, \delta_1  \|T^j y\| \ge \frac{1}{100} \delta_1^2 ,  \text{ since } j < |\log_K \delta_1|
 \end{equation}
 
 By \eqref{eq40},  \eqref{eq41},  and  \eqref{eq42}, the term $(A+B) V_y V_y^* \VV x_0$ in \eqref{eq39} is much more collinear with $y$ than the term 
 $$A V_y V_y^* (x_0 - \VV x_0) \text{ or } A V_y V_y^* x_0.$$
 So, in \eqref{eq39} we must have
 \begin{equation}\label{eq43}
 	\left\| A V_y V_y^* \, (x_0 - \VV x_0)\right\| < 100 \cdot \frac{\gamma'}{\delta_1^2} \cdot (\VE)^{-1/2} \text{ and } \frac{1}{100} \cdot (\VE)^{-1/2} < B < 10 (\VE)^{-1/2}.
 \end{equation}
 We have 
 
$$\frac{1}{10} \cdot (\VE)^{-1/2} < \|y\| < (\VE)^{-1/2} , \frac{1}{10} < \|V_y V_y^* \VV x_0\|< 1 , $$ 
$$\frac{1}{10} \cdot (\VE)^{-1} < \|V_y V_y^* (x_0 - \VV x_0) \|< 10 \cdot (\VE)^{-1} .$$
 
 By \eqref{eq35}--\eqref{eq37} and \eqref{eq28primeprime}--\eqref{eq30primeprime} we have in first order for some linear operator $L$ (determined by the $r_j$'s)
 \begin{equation*}
 \begin{aligned}
 	L(V_y V_y^* \VV x_0) & = \langle ch\VV x_0, x_0 \rangle,\\
 	L(V_y V_y^* [\, \,]^{-2} x_0) & = \langle ch\VV x_0, \VV x_0 \rangle,\\
 	L(\VV y) & = ch \langle \VV x_0, y\rangle
 \end{aligned}
  \end{equation*}
 and, from \eqref{eq39}, 
 
  \begin{equation}\label{eq44}
 	\begin{aligned}
 		| A & \langle ch \VV x_0, x_0 \rangle + B \langle ch \VV x_0, \VV x_0 \rangle - ch \langle \VV x_0, y \rangle | \\
 		& <  \text{ max } \left\{ |\langle ch \VV x_0, x_0\rangle|, |\langle ch \VV x_0, \VV x_0\rangle|, |ch \langle \VV x_0, y\rangle\right\} (\VE)^{98}.
 	\end{aligned}
 \end{equation}

We want $\langle \VV x_0, y \rangle$ after the change to be near to the new $(\VE)^{1/2}$. So, if the change in 
$\langle \VV x_0, y\rangle$ is $-\beta (\VE)^{1/2}$ we want the change in $(\VE)$ to be near to $-2 \beta (\VE)$.

The change in $(\VE)$ in first order is 
$$- ch \langle \VV x_0, x_0 - 2 \VV x_0\rangle .$$

$V_y V_y^* (x_0 - 2 \VV x_0)$ and $y$ are $\VE$-linearly independent so we can make a change $ch \VV x_0$ such that $\VE$ decreases by $-2 \beta(\VE)$. Then, by \eqref{eq39}, 
\begin{equation}\label{eq45}
\langle \VV x_0, x_0\rangle \text{ does not change more than } 10 \beta \VE.
\end{equation}

By \eqref{eq45} and the step in MC just described, it is possible for operators of type 1 to use MC to get down to arbitrarily small $\VE$'s without changing $\langle \VV x_0, x_0\rangle$ more than a total of $10 (\VE)_0$, where $(\VE)_0$ is the start-$(\VE)$ for our MC.

In order to make $\VV x_0 = \VV_n x_0$ converge in norm we need to put is more restrictions on 
$$y_n \mapsto y_n + \sum_{j \ge 0} r_j T^j y_n.$$

In order to do so, we will use that $T$ is of type 1 and that 
$$\left| \langle T^{j_1} y_n, y_n\rangle\right| \ge \delta_1 \|y_n\|^2$$
for some $j_1$ with $1 \le j_1 \le n_0$ and 
$$\left| \langle T^{j_2} y_n, y_n\rangle\right| \ge \delta_2 \|y_n\|^2$$
for some $j_2 > n_0$, with $\delta_1, \delta_2,$ and $n_0$  being independent of $\VE$. Therefore $\delta_1, \delta_2$ have a larger order of magnitude than both $\VE$ and many functions of $\VE$.

We will find a sequence of elements $w_{00}, w_{01}, w_{02},\ldots$ as follows. We first use MC as above without further restrictions than \eqref{eq38}. Then, for a suitable $n_0'$ (we can wait until $(\VE)_n$) is arbitrarily small) we add the condition 

$$\langle ch (\VV_n x_0), w_{00}\rangle = 0$$
 
 for $n > n_0'$. We need to find $w_{00}$ such that the triple 
 $$V_{y_n}V_{y_n}^* (x_0 - \VV x_0), \quad V_{y_n} V_{y_n}^* (w_{00}), \text{ and } y_n$$
are $\delta_2$-independent for a sequence of $n$'s, $n \ge n_0'$ which will allow MC to be used to arbitrarily small $\VE$'s and a small variation of $\VV_n x_0$.

Put 
\begin{equation}\label{eq46}
w_{00} = y'_{n_0'} + \delta_2^{100} \, T^{* j_1} y'_{n_0'} + \delta_2^{100} \, T^{* j_2} y'_{n_0'}.
\end{equation}

Then, since  $\langle T^{* j_1} y_{n}', y_n' \rangle= \langle y_n', T^{j_1} y_n' \rangle$ we get that the triple
$$V_{y_n}V_{y_n}^* (x_0 - \VV x_0), \quad V_{y_n} V_{y_n}^* (w_{00}), \text{ and } y_n$$
is $\delta_2$-linearly independent to arbitrarily small $\VE$'s, since with the condition $\langle ch \VV_n x_0, w_{00}\rangle = 0$ the room for variation of $[\, \,]^{-1}_n x_0$ is less than $\delta_2^{40}$.

We use MC with this restriction until $n = n_1'$ and some small $\VE$, then we put in another restriction, replace $w_{00}$ by $w_{01}$ and 
$$\langle ch [\, \,]^{-1}_n, w_{01}\rangle = 0, \, w_{01} = y'_{n_1'} + \delta_2^{1000} \, T^{* j_1} y'_{n_1'} + \delta_2^{1000}\,  T^{* j_2} y'_{n_1'}$$
similar to \eqref{eq46} but with $\delta_2^{100}$ replaced by $\delta_2^{1000}$.

The triple $V_{y_n} V_{y_n}^* (x_0 - \VV x_0)$, $V_{y_n} V_{y_n}^* (w_{01})$, $y_n$ for $n \ge n_1'$ is $\delta_2$-linearly independent. The room for variation of $[\, \,]^{-1}_n x_0$ is less than $\delta_2^{400}$. By continuing with smaller and smaller $(\VE)$'s, $w_{00}, w_{01}, w_{02}, \ldots$ we have that $V_{y_n} V_{y_n}^* [\,\,]^{-1}_n x_0$ converges to a non-cyclic vector.

It remains to show that operators of type 2 have non-trivial invariant subspaces. Let us assume that \eqref{eq20} holds. Put $$y'_n = y_0' + s_n, \quad \|s_n\| \le D.$$

We will first show that there exists $L_0 > 0$ such that for $L < L_0$, there exists $M(L) = M$ such that if we move $y_0 + L s_n$ within distance $0.3$ of $x_0$ by minimal $\ell_n'$.,
$$\|\ell_{n}' (T) (y_0 + L s_n) - x_0\| \le 0.3$$
then $\|\ell_n'\|_2 \le M$ for all $n$. But, if $L > L_0$ there is no such $M$.
%\textcolor{red}{A PARTIR DE AQUI SIGUE LA PAGINA NUMERADA CON EL 41}.
Let us prove this statement. We first see that $L_0>0$. Since, if 
$$\|\ell'(T) y_0 - x_0\|\le 0.2$$
then
$$\|\ell'(T) (y_0' + L s_n) - x_0\| < 0.3 \text{  if  } L < \frac{0.1}{D \, \|\ell' (T)\|_2}.$$
Assume then that, for every $L$, there exists $M(L)$ such that for every $n$, there exists $\ell'_n(T),$ $\|\ell'_n(T)\|_2 \le M(L)$ such that
$$\|\ell'(T) (y_0' + L s_n) - x_0\| \le 0.3.$$
With $\|T\|_{\text{op}} = \frac{1}{K}$, $K= 10^{20},$ the terms of degree larger than 
$$10 \log_{K} (M(L) + D \, L \, + \|y\|) = M'$$
in $\ell_n'(T)$ will contribute less than $1/10^{20}$ to 
$$\sum_{j \ge 0} a_j T^j (y_0' + L s_n) = \ell_n' (T) (y_n' + L s_n).$$
For terms with degree less than or equal to $M'$, by \eqref{eq20}, 
$$\left| \langle \sum_{j=0}^{M'} a_j T^j y_0', T^{*m} q(T^*) y_0'\rangle + \frac{1}{L} 
\langle \sum_{j=0}^{M'} a_j T^j (L s_n),  T^{*m} q(T^*) s_n\rangle \right| \longrightarrow 0 \text{ as } n \rightarrow \infty.$$
With $T^{*m} {q} (T^*)y_0'$ near $x_0$ this shows that 
$$\|\ell_n(T) (y_0 + L s_n)\| > 1 + \frac{1}{5}.$$
This contradiction shows the statement above.

We now let $\delta >0$, $L < L_0$ such that $|L - L_0| < \delta/L$, put $M = M(L)$ and
$$\|\ell_n(T) (y_0 + L s_n) - x_0\| \le 0.3 \text{ and } \|\ell_n(T)\|_2 \le M \text{ for all } n.$$
Then, 
$$\|\ell_n(T) (y_0 + L_0 s_n) - x_0\| = \|\ell_n(T) (y_0 + L s_n) + (L_0-L) s_n - x_0\|.$$

Since $\|\ell_n(T) L s_n\| < 2,$  $\|(L_0 - L) \ell_n(T) s_n\| < 2 \delta$ and, so, 
$$\|\ell_n(T) (y_0 + L_0 s_n) - x_0 \| < 0.3 + 2 \delta \, \| \ell_n(T)\|_2 \le M.$$
 
Now, if $M(L_0) = \infty$ then, for some $n$, 
$$\|\ell(T) (y_0' - L_0 s_n) - x_0\| \le 0.3$$
may require an arbitrarily large $\|\ell (T)\|_2$.

So, for moving $y_0 + L_0 s_n$ by minimal $\ell(T)$'s between the distances $0.3 + \delta$ and $0.3$ to $x_0$, may require arbitrarily small $(\VE)$'s.

If $M(L_0) < \infty$, then for every $n$, $\|\ell_n(T)\|_2 \le M(L_0)$ for having $\|\ell_n(T) (y_0 + L_0 s_n) - x_0\| \le 0.3$.
%\begin{align}\label{eq47}
%\hspace{-.5cm}\begin{split}	
%& \text{If the} (\VE)'s \text{ here were bounded from below} (\VE)_n > 0 \text{ we could find,}\\
%& \text{for some } L > L_0, \text{ move } y_0 + L s_n \text{ within } 0.3 - \delta',  \delta' > 0, \\
%& \text{with a uniform bound on the } \| \ell_n(T)\|'s, \text{ which would be a contradiction.}\\ 
%&\text{So, also here we get arbitrarily small } (\VE)'s.
%\end{split}
%\end{align}

If the $(\VE)$'s here were bounded from below $(\VE)_n > 0$ we could find, for some $L > L_0$, move $y_0 + L s_n$ within $0.3 - \delta'$, $\delta' > 0$, with a uniform bound on the $\| \ell_n(T)\|$'s, which would be a contradiction. So, also here we get arbitrarily small $(\VE)$'s.\hfill (47)

Thus, we will find a sequence  $y_n' = y_0' + L_0 s_n$ and $\delta_n > 0$, $\delta_n \rightarrow 0$, such that we have: With $\|\ell'_n\|_2$ being inf $\|\ell'\|_2$ for 
$$\|\ell'(T) (y_0' + L_0 s_n) - x_0\| \le 0.3 \text{  (or } \le 0.3 + \delta_n \text{)}$$
we put 
$$\ell'(T) (y_0' + L_0 s_n) = y'_n.$$
We have, for each $n$, $\langle y'_{n_1}, x_0 - y'_{n_1}\rangle = (\VE)_n$ with $(\VE)_n$ tending to $0$ and 
$$\| x_0 - y'_{n_1}\| = 0.3 \text{  (or } 0.3 + \delta_n \text{)}.$$
We use Lemma 2 in order to find minimal $\|\ell'_{n_1}\|_2$ such that 
$$\|\ell_{n_1}'(T) y'_{n_1} - x_0\| \le 0.3 \text{  (or } \le 0.3 + \delta_n \text{)}$$
with dominating $|a_0|$ in $\ell'_{n_1}(T) = \sum_{j \ge 0} a_j T^j.$ We apply MC with
$$y'_{ni+1} = y'_{n i} + \sum_{j \ge 0} r_j T^j y'_{n i}$$
and continue MC, $i= 1, 2, 3, \ldots$ until either, for some $m$,
\begin{itemize}
	\item[1.)] $\text{sup}_{j \ge m} |\langle T^j y'_{n i}, y'_{n i}\rangle | \ge \gamma^{1/10} \, (\VE)_n$ for each $i$. Then we continue MC as for operators of type 1, $n \rightarrow \infty$, and we get a non-cyclic vector.
	\item[2.)] For some $i$, 
	$$|\langle T^j y'_{n i}, y'_{n i}\rangle | < \gamma^{1/10} \, (\VE)_n$$
	for all $j \ge m$.
\end{itemize}
Then we start over with a new $y'_n$ (call it $y'_{n+1}$) 
$$y'_{n+1} = y'_0 + L_0 s_{n+1},$$
and a new $(\VE)_{n+1}$ which is much smaller than $(\VE)_n$.

If, for every $n$, when we have started over, 2.) happens for some $i = i(n)$, which may vary with $n$, then for $y'_{n i}$ we have, for $j \ge m$,
\begin{align*}
|\langle T^j y'_{n i}, x_0\rangle| & = | \langle T^j y'_{n i}, y'_{n i} + x_0 - y'_{n i}\rangle | \\
& \le | \langle T^{j} y'_{n i} , y'_{n i}\rangle| + |\langle T^j y'_{n i} , x_0 - y'_{n i} \rangle |\\
& =  | \langle T^{j} y'_{n i} , y'_{n i}\rangle| \\
& + |\langle T^j \ell'_{n i}(T) y'_{n i} + y'_{n i} - \ell'_{n i}(T) y'_{n i}, x_0 - \ell'_{n i}(T) y'_{n i} + \ell'_{n i}(T) y'_{n i} - y'_{n i}  \rangle|\\
& \le |\langle T^j y'_{n i}, y'_{n i}\rangle| + |\langle T^j \ell'_{n i}(T) y'_{n i}, x_0 - \ell'_{n i} (T) y'_{n i}  \rangle|\\ 
& +  \|T^j (y'_{n i} - \ell'_{n i} (T) y'_{n i})\| + \|\ell_{n i} (T) y'_{n i} - y'_{n i}\| \longrightarrow 0
\end{align*}
for every $j \ge m$ as $n \rightarrow \infty$ and, for $w$-$\displaystyle\lim_{n \rightarrow \infty} y'_{n i} = y'_\infty$ we get $T^m y'_\infty$ non-cyclic, since $\langle T^j y'_\infty, x_0\rangle = 0$ for $j \ge m$.
\bigskip

Notice, as a final remark, that even if $y'_n$ tends weakly to $y_0$ as $n \rightarrow \infty$, we should not expect $y'_\infty $ and $y_0$ to be the same.\\

%%%%%%%%%%%%%%
%%%%%%%%%%%%%%
%%%%%%%%%%%%%%
%%%%%%%%%%%%%%
%%%%%%%%%%%%%%

\noindent \textbf{Summary of Construction.}   We will find a sequence $(y'_n)$ in $H$, an $x_0 \in H$ with $\|x_0\| =1$, and $\varepsilon \in (0,1)$. For every  $n$ we will find an analytic expression
$$\ell_n (T) = \sum_{j \ge 0} a_j T^j,$$
such that 
$$\|x_0 - \ell_n (T) y_n)\| \le \varepsilon$$
and such that 
$$\| \ell_n(T) \|_2^2 = \sum_{j\ge 0} |a_j|^2$$
is minimal for this. We will, in the Main Construction have
$$y_{n+1} = y_n + \sum_{j \ge 0} r_j T^j y_n$$
but sometimes we define $y_{n+1}$ from $y_n$ in another way.

Most of the paper, Lemma 1--Lemma 5 and \eqref{eq01}--\eqref{eq27} are concerned with the task of finding suitable $x_0$, $y_0$, and $\ell_0$ such that the Main Construction can start and continue.

We put $\displaystyle \VE = \langle \ell(T) y', x_0 - \ell (T) y' \rangle$. By \eqref{eq02}--\eqref{eq08}, $\VE \ge 0$, and if $\VE =0$ we have already an invariant subspace.

We introduce a representation of the minimal $\ell$ ($ = \ell$ of minimal $\|\cdot\|_2$-norm) by introducing $V_y$ and $V_y^*$, and using $V_y V_y^*$, \eqref{eq02}--\eqref{eq13}.

To start the Main Construction we need $\displaystyle \ell'(T) = \sum_{j \ge 0} a_j T^j$ to have a dominating $|a_0|$, and we use Lemma 2 to achieve that. In Lemma 2 we start with a general $\ell (T) y' = y'_1$, $\|\ell(T)\|_2$ minimal to have
$$\|x_0 - \ell(T) y'\| = \|x_0 - y_1'\|.$$
And Lemma 2 gives us, if $\VE$ is sufficiently small that, if $\displaystyle \|x_0 - \sum_{j \ge 0} a_j T^j y_1'\| = \|x_0 - y_1'\|$ with minimal $\displaystyle \sum_{j \ge 0} |a_j|^2$, then $|a_0|$ is dominating in $\displaystyle \sum_{j \ge 0} a_j T^j$. This is at the cost of possibly increasing $\VE$, but that will be compensated for, when the Main Construction is used later.

We also see, from \eqref{eq14} and \eqref{eq17}, that: If 
$$\|x_0 - \ell' (T) y' \| = \|x_0 - V_y V_y^* \VV x_0\|$$ 
and $\|y'\|$ is of order of magnitude $1$ and $|a_0'|$ is dominating in  $\displaystyle \ell'(T) = \sum_{j \ge 0} a'_j T^j y'$, then $\|y\|$ is of order of magnitude $(\VE)^{-1/2}$. To achieve, that we can use Lemma 2, we need a sufficiently small $\VE$ to start with. In order to get that, we use Lemma 1, \eqref{eq19}--\eqref{eq21} and $u_0$ and $u_1$ from the definition of operators of type 1 and type 2. When we make the definition of type 1 and type 2 operators we also make the choice of $u_0$ and $u_1$, $y_0 = u_0$, and $\displaystyle x_0 = \frac{u_0}{2} + \frac{\sqrt{3}}{2} u_1.$

So, by using Lemmas 1 -- 5, and \eqref{eq11}--\eqref{eq24} we have either an invariant subspace or a situation with 
$$\|x_0 - \ell'(T) y' \le \varepsilon\|, \, \ell' (T) = \sum_{j \ge 0} a_j T^j \, \text{ with } \sum_{j \ge 1} |a_j|^2 \le \VE \cdot \gamma (\VE)$$
with $\VE$ arbitrarily small. Then, either we an estimate from below $\displaystyle \sum_{j \ge 1} |a_j|^2 \ge f(\VE)$ where $f(\VE)$ can be small, say $f(\VE) = (\VE)^{100}$, and we have what we need, or  we do not have such an estimate. Then, by using Lemmas 1--5 and \eqref{eq24}--\eqref{eq27} with Case I and Case II we either get an invariant subspace or an estimate \eqref{eq25} from below of $\displaystyle \sum_{j \ge 1} |a_j|^2 \ge (\VE)^{15}$, and the estimate from above is automatic.

For  $\displaystyle \sum_{j \ge 1} |b_j|^2$ we can use the same argument to get an estimate from below.

So, with this preparations, we  can start the Main Construction. By \eqref{eq28primeprime}--\eqref{eq30primeprime}, \eqref{eq31}--\eqref{eq34} we get the needed product estimates for analytic function of $(\VE)$, having some concentration at degrees not higher than some function of $(\VE)$ (it follows from \cite{BeauzamyEnflo}, but it can also be done in a completely elementary way).  We use that $y$ and $V_y V_y^* \VV x_0$ are $\gamma (\VE)$-close to being collinear but $V_y V_y^* (y)$ and $y$ are $\delta_2$-far from being collinear, since $T$ is of type 1. By \eqref{eq35}--\eqref{eq45} we can get down to  arbitrarily small $(\VE)$'s by the Main Construction.

To get $\ell_n(T) y'_n$ to converge we put in more restrictions of changes when $y_n \rightarrow y_{n+1}$. By \eqref{eq46} we let the changes of $\VV x_0$ be orthogonal also to $w_{00}$ and then $w_{01}, w_{02}$, etc. Since the triplet 
$$V_{y_n} V_{y_n}^* (x_0 - \VV x_0), V_{y_n} V_{y_n}^* (w_{0j}), \text{ and } y_n$$
will be $\delta_2$-linearly independent for every $j$ and, with increasing $j$, $\ell_n (T) y_n$ can vary less and less so that,  with $\VE \rightarrow 0$, $\ell_n(T) y'_n$ converges for a non-cyclic vector.

For operators of type 2 we use \eqref{eq20} and the argument from (47) in order to produce arbitrarily small $(\VE)$'s. Then we use Lemma 2 and \eqref{eq24}--\eqref{eq27} to start the Main Construction, starting with a $y_n$, $(\VE)_n$, and 
$$\text{sup} |\langle T^j y_n, y_n\rangle | > \delta \|y_n\|^2.$$
When we do the Main Construction we do not control $\text{sup} |\langle T^j y_n, y_n\rangle |$ so the Main Construction gives us vectors $y_{n1}, y_{n2}, \ldots$ 

Case I: If, for all $i$,  
$$\text{sup}\left| \langle T^j y_{ni}, y_{ni}\rangle\right| \ge \left( \gamma ((\VE)_n)\right)^{1/10} \|y_{ni}\|^2$$
we have the same situation as for operators of type 1 and we continue to a non-cyclic vector.
%%%%%%%%%%%%%%

Case II: If, for some $i$, 
$$\text{sup}\left| \langle T^j y_{ni}, y_{ni}\rangle\right| < \left( \gamma ((\VE)_n)\right)^{1/10} \|y_{ni}\|^2$$
we interrupt the process and start over with a much smaller $(\VE)_{n+1}$. 
%%%%%%%%%%%%%%

If Case I happens for some $y_n$ and the sequence $y_{ni}$ we get an invariant subspace. If, for every $n$, the Case II happens for some $y_{ni}$, $i = i(n)$, then for a weakly convergent subsequence in $n$ of $y'_{ni}$ will have $\displaystyle T^m (\lim_{n \rightarrow \infty} y'_{ni})$ non-cyclic for some $m$.\\

%%%%%%%%%%%%%%
%%%%%%%%%%%%%%
\bigskip
%%%%%%%%%%%%%%
%%%%%%%%%%%%%%
\noindent \textbf{Acknowledgments.} I would like to sincerely thank professors  \linebreak Juan~B.~Seoane~Sep\'ulveda and Gustavo A. Mu\~{n}oz Fern\'andez  for invaluable discussions and help to get my manuscript in order. Sincere thanks also to professors \linebreak A. Aleman, E. Rydhe and A. Bergman for going over the earlier version of the paper with me and providing  very valuable comments. Thanks, as well, to professors \linebreak B. Karltag, A. Sodin, A. Alpeev and G. Schechtman for their work and comments on this earlier version. For their long-term support for this whole project I sincerely thank professors W. B. Johnson, M. S. Moslehian and J. Virtanen. Last, but not least, I want to thank my wife Lena for her long and unfailing support.

%%%%%%%%%%%%%%
%%%%%%%%%%%%%%
% REFERENCES
%%%%%%%%%%%%%%
%%%%%%%%%%%%%%
%%%%%%%%%%%%%%
%%%%%%%%%%%%%%

\end{document}